\documentclass[12pt]{article}
\usepackage[colorlinks,
            linkcolor=blue,
            anchorcolor=green,
            citecolor=blue
            ]{hyperref}
\usepackage{mathrsfs,amssymb,amsfonts,soul}
\usepackage{amsmath,amssymb,latexsym,color}
\usepackage[mathscr]{eucal}
\usepackage[numbers,sort&compress]{natbib}
\usepackage{graphicx}
\usepackage{lipsum}
\usepackage{tikz,tikz-cd}
\usepackage{caption}

\usetikzlibrary{matrix, arrows, automata, positioning, calc}
\def\thefootnote{\fnsymbol{footnote}}
\newtheorem{thm}{Theorem}[section]

\newtheorem{prop}[thm]{Proposition}
\newtheorem{lemma}[thm]{Lemma}
\newtheorem{cor}[thm]{Corollary}
\newtheorem{example}[thm]{Example}

\newtheorem{problem}[thm]{Problem}
\newtheorem{ques}[thm]{Question}
\newtheorem{remark}[thm]{Remark}
\newtheorem{Notation}[thm]{Notation}

\newcommand{\proof}{{\noindent{\it Proof.\quad}}}
\newcommand{\qed}{$\hfill\Box\medskip$} % <<<<<<<<<<<< I PUT DOLLARS HERE, IN ORDER TO WRITE \qed INSTEAD OF $\qed$ IN THE TEX.
\renewcommand{\thefootnote}{\arabic{footnote}}

\newcommand{\Cyc}{{\rm Cyc}}

\usepackage{CJK}

\newcommand{\sime}{\stackrel{e}{\sim}}
\newcommand{\equive}{\stackrel{e}{\equiv}}
\usepackage[english]{babel}
\usepackage{csquotes}
\usepackage{enumitem}
\setlist[enumerate,1]{label={\rm{(\alph*)}}}
\begin{document}
\begin{CJK*}{GBK}{song}
\renewcommand{\abovewithdelims}[2]{
\genfrac{[}{]}{0pt}{}{#1}{#2}}
%%%%%%%%%%%%%%%%%%%%%%%%%%%%%%%%%%%%%%%%%%%%%%%%%%%%%%%%%%%%%%%%%%%%%%%%%%%%%%%%%%%%%%%%
%%%%%%%%%%%%%%%%%%%%%%%%%%%%%%%%%%%%%%%%%%%%%%%%%%%%%%%%%%%%%%%%%%%%%%%%%%%%%%%%%%%%%%%%

\title{\bf Forbidden subgraphs in enhanced power graphs of finite groups}

\author{{Xuanlong Ma$^{1}$, Samir Zahirovi\'c$^{2,}$\footnote{Corresponding author}~,
Yubo Lv$^{3}$, and Yanhong She$^{1}$}
\\
{\small\em $^1$School of Science, Xi'an Shiyou University, Xi'an,  710065, China}
\\
{\small\em $^2$Department of Mathematics and Informatics, Faculty of Sciences,} 
\\
{\small\em University of Novi Sad, Dositej Obradovi\'c Square 3, Novi Sad, 21000, Serbia}
\\
{\small\em $^3$School of Mathematical Sciences,} 
\\
{\small\em Guizhou Normal University, Guiyang, 550001, China}
}

 \date{}
 \maketitle
\newcommand\blfootnote[1]{%
\begingroup
\renewcommand\thefootnote{}\footnote{#1}%
\addtocounter{footnote}{-1}%
\endgroup
}
\begin{abstract}
The enhanced power graph of a group is the simple graph whose vertex set is consisted of all elements of the group, and whose any pair of vertices are adjacent if they generate a cyclic subgroup. In this paper, we classify all finite groups whose enhanced  power graphs are split and threshold.
We also classify all finite nilpotent groups whose enhanced  power graphs are chordal graphs and cographs. Finally, we give some families of non-nilpotent groups whose enhanced power graphs are chordal graphs and cographs. These results partly answer a question posed by Peter J. Cameron.


\end{abstract}

\blfootnote{{\noindent E-mail addresses: xuanlma@mail.bnu.edu.cn (X. Ma), samir.zahirovic@dmi.uns.ac.rs\\ (S. Zahirovi\' c), lvyubo341281@163.com (Y. Lv), yanhongshe@xsyu.edu.cn (Y. She)}}

\section{Introduction}

The {\em undirected power graph} $\mathcal{P}(G)$ of a group $G$ has vertex set $G$, and two distinct elements of $G$ are adjacent in the undirected power graph if one is a power of the other. The concept of the directed power graph was first introduced by Kelarev and Quinn \cite{KQ00}, and the undirected power graph was first studied by Chakrabarty et al. \cite{CGS}.
In recent years, the interest among researchers in the undirected power graph has been growing; see, for example, \cite{B3,Cam,CGh,FMW,MWW,PaG}. Also, see \cite{AKC, Came2021, KSCC} for a survey of results and open problems on the undirected power graph.

Let $G$ be a finite group. The {\em enhanced power graph} of $G$, denoted by $\mathcal{P}_e(G)$, is the simple graph whose vertex set is the set of all elements of $G$, and two distinct vertices $x,y$ are adjacent if and only if $\langle x,y\rangle$ is a cyclic subgroup of $G$. The term enhanced power graph was introduced by Aalipour et al. \cite{acam}; in their paper, they studied the enhanced power graph in order to measure how close the undirected power graph is to the commuting graph. It is important to mention that this notion was studied even earlier in \cite{afkami, ma-wei-zhong}. It is worth noting that in  \cite{afkami, ma-wei-zhong}, as well as in some newer papers (e.g. \cite{CLSTU,CLSTU2}), the authors used the term {\em cyclic graph} instead. Furthermore, Abdolahi et al. \cite{abdolahi} studied the noncyclic graph, i.e. the complement of the enhanced power graph, even earlier.

Recently, the enhanced power graphs has received considerable attention. Aalipour et al. \cite{acam} showed that every maximal clique of $\mathcal{P}_e(G)$ is a maximal locally cyclic subgroups of $G$.
In 2020, Zahirovi\'c, Bo\v{s}njak, and Madar\'asz \cite{Za} showed that any isomorphism between undirected power graph of finite groups is an isomorphism between enhanced power graphs of these groups. Panda, Dalal, and Kumar \cite{Pa} studied the minimum degree, independence number, and matching number of enhanced power graphs of finite groups. In \cite{MS}, the authors  obtained an explicit formula for the metric dimension of the enhanced power graph of a finite group. Bera, Dey, and Mukherjee \cite{Bera} completely characterized all abelian groups such that their proper enhanced power graphs are connected, where the proper enhanced power graph of $G$ is the graph obtained by deleting the identity element from $\mathcal{P}_e(G)$. Dalal and Kumar \cite{Da} studied the enhanced power graphs of some special groups. Further, Cameron and Kuzma \cite{CamK} defined a graph whose vertex set is a finite group $G$, whose edge set is contained in that of the commuting graph of $G$ and contains the enhanced power graph of $G$; they called it the deep commuting graph. In \cite{CamK}, they obtained conditions for the deep commuting graph to be equal to either of the enhanced power graph and the commuting graph.

A number of important graph classes, including cographs, chordal graphs, split graphs, and threshold graphs, can be defined either structurally or in terms of forbidden induced subgraphs. Forbidden subgraphs of undirected power graphs of groups have been studied by Doostabadi et al. \cite{Doo2014} and Cameron et al. \cite{Cam2020}.

In \cite{Came2021}, Cameron described various graphs whose vertex set is the universe of a group and whose edges reflect group structure in some way, and he proposed the following question.

\begin{ques}\label{ques1}{\rm (\cite[Question 14]{Came2021})}
For which finite groups is the enhanced power graph a perfect graph, a cograph, a chordal graph, a split graph, or a threshold graph?
\end{ques}

Recently, Bo\v{s}njak, Madar\'asz, and Zahirovi\'c \cite{Bos,Za} studied  perfectness of the enhanced power graph of a group. Motivated by Question \ref{ques1}, in this paper, we study forbidden subgraphs of enhanced power graphs of finite groups. This way,
we classify all finite groups whose enhanced  power graphs are split graphs and threshold graphs.
We also classify all finite nilpotent groups whose enhanced  power graphs are chordal graphs and cographs.
Finally, we give some families of non-nilpotent groups whose enhanced power graphs are chordal graphs and cographs. Our results partly answer Question \ref{ques1}.

\section{Preliminaries}

In this section, we briefly recall some notation, terminology, and basic results which we need in the remainder of this paper.

All groups considered in this paper are finite.
We denote groups by capital letters, usually $G$, and its identity element is usually denoted by $e$. The {\em order} of a group element $g$ is denoted by $o(g)$, the set of orders of elements of a group $G$ is denoted by $\pi_e(G)$. An element of order $2$ is called an {\em involution}.
By $\mathcal{M}_G$, we denote the set of all maximal cyclic
subgroups of $G$. Notice that $|\mathcal{M}_G|=1$ if and only if $G$ is cyclic.
The cyclic group of order $n$ is denoted by $\mathbb{Z}_n$. For a prime number $p$, a finite elementary $p$-group is an abelian group of exponent $p$, and it is isomorphic to $\mathbb{Z}_p^n$ for some $n\in\mathbb N$.

We remind the reader that a finite nilpotent group is the direct product of its Sylow subgroups. Also, it is well known that both $p$-groups and abelian groups are nilpotent.

By $D_{2n}$, we denote the dihedral group of order $2n$. Let us recall that $D_{2n}$ consists of $n$ reflection, $n-1$ rotations, and the identity mapping. Together with the identity mapping, all rotations make up a cyclic subgroup of order $n$, and each reflection makes up a maximal cyclic subgroup of order $2$ of $D_{2n}$. Dihedral group $D_{2n}$ has presentation
\begin{equation}\label{d2n}
D_{2n}=\langle a,b: a^n=b^2=e, bab=a^{-1}\rangle.
\end{equation}
Notice that, for $a$ and $b$ from the above presentation, $\langle a\rangle$ is the cyclic subgroup containing all rotations, and $a^ib$ is a reflection for every $i$. Therefore,
\begin{equation}\label{d2n-2}
\mathcal{M}_{D_{2n}}=\{\langle a\rangle,\langle b\rangle,\langle ab\rangle,
\langle a^2b\rangle,\ldots,\langle a^{n-1}b\rangle\}.
\end{equation}

The generalized quaternion group $Q_{4m}$ of order $4m$, $m\ge 2$, is the group defined by presentation
\begin{equation}\label{q4m}
Q_{4m}=\langle x,y: x^m=y^2, x^{2m}=e, y^{-1}xy=x^{-1}\rangle.
\end{equation}
Let us remark that, for $x$, $y$ and $m$ from the above presentation, $x^m$ is the unique element of $Q_{4m}$ of order $2$. Also, $\langle x\rangle$ is its maximal cyclic subgroup of order $2m$, and $\langle x^iy\rangle$ is a maximal cyclic subgroup of order $4$ for each $i$, $0\le i<m$. Therefore,
\begin{equation}\label{q4m-2}
\mathcal{M}_{Q_{4m}}=
\{\langle x\rangle,\langle xy\rangle,
\langle x^2y\rangle,\ldots,\langle x^my\rangle\}.
\end{equation}

All graphs considered in this paper are finite, undirected, and with no loops and no multiple edges.
Let $\Gamma$ be a graph.
The vertex set of $\Gamma$ is denoted by $V(\Gamma)$,
and the edge set by $E(\Gamma)$.
A graph is {\em $\Gamma$-free} if it has no induced subgraph isomorphic to $\Gamma$.
A graph is a {\em star} if it is a tree with $n$ vertices in which one vertex has degree $n-1$ and every other
vertex has degree $1$.
We denote by $P_n$, $C_n$, and $K_n$ the path with $n$ vertices, the cycle of length $n$, and the complete graph of order $n$, respectively. Note that the length of a path with $n$ vertices is $n-1$. Also,
we use $2K_2$ to denote the disjoint union of two copies of $K_2$. If two distinct vertices $x$
and $y$ of $\Gamma$ are adjacent, then we denote this by
$x \sim_{\Gamma} y$, or shortly by $x \sim y$.
In particular, a path with $n$ vertices is usually denoted by $x_1\sim x_2\sim\cdots \sim x_n$, and a cycle of length $n$ is denoted by $x_1\sim x_2\sim\cdots \sim x_n
\sim x_1$, where $x_i\in V(\Gamma)$.

By $x\sime y$, we denote that the elements $x$ and $y$ of a group are adjacent in the enhanced power graph. The {\it closed neighborhood} of a vertex $x$ in a graph $\Gamma$ is the set $\{x\}\cup\{v\in V(\Gamma): v\sim_\Gamma x\}$. If elements $x$ and $y$ of a group $G$ have the same closed neighborhood in $\mathcal P_e(G)$, we write $x\equive y$; in that case, we also say that $x$ and $y$ are {\it closed twins} in $\mathcal P_e(G)$.

In this paper, we study several properties of the enhanced power graph of a group which are characterizable using \enquote{forbidden subgraphs}. We introduce now their definitions, as well as their characterizations by forbidden subgraphs. A {\em split graph} is a graph whose vertex set is the disjoint union of two subsets $A$ and $B$ such that $A$ induces a complete graph and $B$ induces an empty graph. It is known that a graph is split if and only if it does not contain $C_4$, $C_5$, nor $2K_2$ as an induced subgraph. A graph is called a {\em threshold graph} if it can be constructed by addition of an isolated vertex and addition of a dominating vertex, i.e. a vertex adjacent to all vertices. Equivalently, a graph is threshold if and only if it has no induced subgraph isomorphic to $P_4$, $C_4$, or $2K_2$. A graph is called {\em chordal} if it contains no induced cycles of length greater than $3$. In other words, a chordal graph is a graph in which every cycle of length at least $4$ has a chord, i.e. an edge incident to a pair of nonconsecutive vertices of the cycle. Notice that, if a graph is $C_4$-free and $P_4$-free, then it is chordal. Note, however, that this sufficient condition is not a necessary one as well. A graph is called a {\em cograph} if it can be constructed from isolated vertices by disjoint union and complementation. It is known that a graph is a cograph if and only if it is $P_4$-free.

The relations between the above-introduced families of graphs is presented by Diagram \ref{hase-diag-1}; it is well known that cographs and split, threshold, and chordal graphs are all perfect graphs.  Note that perfectness of a graph is also a property which, by the Strong Perfect Graph Theorem, is characterizable by forbidden subgraphs.

\renewcommand{\figurename}{Diagram}

\begin{figure}[h!]
\centering
\begin{tikzpicture}
  \node (a) at (0,0) {threshold graphs};
  \node (b) at (-2,1.5) {cographs};
  \node (c) at (0,3) {perfect graphs};
  \node (d) at (2,1) {split graphs};
  \node (e) at (2,2) {chordal graphs};
%  \path[-stealth, auto] (a) edge[draw=none]  node [sloped, auto=false, allow upside down] {$\subset$} (b)
%                                     edge[draw=none]  node [sloped, auto=false, allow upside down] {$\subset$} (d)
%                                (b) edge[draw=none]  node [sloped, auto=false, allow upside down] {$\subset$} (c)
%                                (d) edge[draw=none]  node [sloped, auto=false, allow upside down] {$\subset$} (e)
%                                (e) edge[draw=none]  node [sloped, auto=false, allow upside down] {$\subset$} (c);
  \draw (a) -- (b) -- (c) -- (e) -- (d) -- (a);
%\path[-] (a) edge[bend left=15]  (b)
%                    (b)    edge [bend left=15]   (c)
%                     (c)     edge [bend left=15] (e)
%                    (e) edge [bend left=15]   (d)
%                      (a)     edge [bend right=15]   (d);
\end{tikzpicture}
\caption{}  \label{hase-diag-1}
\end{figure}

Notice that the cyclic subgroup generated by a pair of elements $x$ and $y$ of a group $G$ is included in every subgroup of $G$ containing $x$ and $y$. Therefore, it is clear that, for any subgroup $H$ of $G$, $\mathcal{P}_{e}(H)$ is an induced subgraph of $\mathcal{P}_{e}(G)$.
The proof of the following result is straightforward.

\begin{lemma}\label{nnpg-lem-new}
Let $G$ be a group and $H$ its subgroup. If $\mathcal{P}_{e}(G)$ is split, threshold, chordal, chordal, or a cograph, then $\mathcal P_e(H)$ is also split, threshold, chordal, chordal, or a cograph, respectively.
\end{lemma}

We introduce now the following important lemma proven by Aalipour et al.\cite{acam}.

\begin{lemma}{\rm (\cite[Lemma 33]{acam})}\label{complete}
Let $G$ be a finite group, and let $C$ be a maximal clique of $\mathcal{P}_{e}(G)$. Then $C$ is a maximal cyclic subgroup of $G$.
\end{lemma}

\section{Split graphs and threshold graphs}\label{}

In this section, we classify all groups whose enhanced power graphs are split and threshold graphs. Before proving Theorem \ref{ep-thm-split}, we shall prove the following two lemmas, which are required for their proofs.

\begin{lemma}\label{ep-sg-lem1}
Let $G$ be a group with $2K_2$-free enhanced power graph. Then $G$ has at most one maximal cyclic subgroup of order greater than $2$.
\end{lemma}

\proof
If $G$ is a cyclic group or has exponent $2$, then the assertion of the lemma obviously holds. So suppose further that $G$ is a noncyclic group.

Suppose, for contradiction, that $G$ has two distinct maximal cyclic subgroups, say $\langle x\rangle$ and $\langle y\rangle$, of orders greater than $2$. Since $\langle x\rangle,\langle y\rangle<\langle x,y\rangle$, and since $\langle x\rangle$ and $\langle y\rangle$ are maximal cyclic subgroups, the group $\langle x,y\rangle$ is not cyclic. Thus, we get $x\not\sime y$, $x\not\sime y^{-1}$, $x^{-1}\not\sime y$, and $x^{-1}\not\sime y^{-1}$. Therefore, the subgraph of $\mathcal{P}_{e}(G)$ induced by $\{x,x^{-1},y,y^{-1}\}$ is isomorphic to $2K_2$, which is in contradiction with the assumption that $\mathcal P_e(G)$ is $2K_2$-free. This proves the lemma.
\qed

\begin{lemma}\label{ep-sg-lem2}
Let $G$ be a group. Then
$\mathcal{P}_{e}(G)$ is $2K_2$-free if and only if $G$ is a cyclic group, a dihedral group, or an elementary abelian $2$-group.
\end{lemma}

\proof
Since $\mathcal{P}_{e}(\mathbb Z_2^n)$ is a star for every $n\geq 2$, $\mathcal{P}_{e}(\mathbb Z_2^n)$ is $2K_2$-free. Also, since the enhanced power graph of a cyclic group is complete, it is $2K_2$-free too. Now, it only remains to show that $\mathcal P_e(D_{2n})$ is also $2K_2$-free. Suppose that some four-element set $D$ induces a subgraph of $\mathcal P_e(D_{2n})$ isomorphic to $2K_2$. Notice that $e\not\in D$ since $e$ is adjacent to all elements of the group. Furthermore, $D$ contains no reflection, since every reflection is adjacent to no nonidentity element. Therefore, $D$ contains only rotations, but that implies that $D$ induces a complete subgraph of $\mathcal P_e(D_{2n})$, which is a contradiction. This proves one implication of the lemma.

Suppose now that $\mathcal{P}_{e}(G)$ is $2K_2$-free. It suffices to prove that if $G$ is noncyclic, then $G$ is either a dihedral group or an elementary abelian $2$-group. Now, let us assume that $G$ is not cyclic. If every element of $G$ is an involution, then $G$ is an elementary abelian $2$-group, as desired. Therefore, we may assume that $G$ has an element of order greater than $2$. Then by Lemma \ref{ep-sg-lem1}, $G$ has a unique maximal cyclic subgroup $\langle a\rangle$ of order is at least $3$. Note that $\langle a\rangle$ contains all elements of $G$ of order greater than $2$. Now, let $b\in G\setminus \langle a\rangle$. Then $ab\notin \langle a\rangle$. It follow that $o(b)=o(ab)=2$ since $ab\not\in\langle a\rangle$, and so
\begin{equation}\label{ep-eq-1}
bab=b(ab)=bb^{-1}a^{-1}=a^{-1}.
\end{equation}
We deduce that $G$ is nonabelian. Moreover, \eqref{ep-eq-1} also implies that $\langle a\rangle\unlhd G$, and hence $\langle a\rangle\langle b\rangle$ is a subgroup of $G$.

In order to prove that $\langle a\rangle\langle b\rangle=G$, suppose that $g\in G\setminus \langle a\rangle$, and let us prove that $g\in \langle a\rangle\langle b\rangle$. It is sufficent to prove that $gb\in \langle a\rangle$. Suppose, for  contradiction, that $gb\notin \langle a\rangle$. Then $o(g)=o(gb)=2$ and $gbg=gg^{-1}b^{-1}=b^{-1}$, which implies that $gb=bg$. Similarly, we have that $xax=a^{-1}$ for every $x\in G\setminus\langle a\rangle$. Also, since $gba\notin \langle a\rangle$, $gba$ is an involution. As a result, we have
$$a^{-1}=a^{-1}(gb)(gb)=\big((gb)a\big)^{-1}(gb)=(gba)gb=g(bab)g=ga^{-1}g=a,$$
which is a contradiction as $o(a)\ge 3$.

We conclude that $g\in \langle a\rangle\langle b\rangle$, which implies that $G=\langle a\rangle\langle b\rangle=\langle a,b\rangle$. Thus, by \eqref{ep-eq-1} and the presentation of the dihedral group, it follows that $G\cong D_{2n}$ where $n=o(a)$. This proves the lemma.
\qed

Now, we are ready to prove the main theorem of the section.

\begin{thm}\label{ep-thm-split}
The following are equivalent for a finite group $G$:
\begin{enumerate}
\item $\mathcal{P}_{e}(G)$ is split;
\item $\mathcal{P}_{e}(G)$ is threshold;
\item $\mathcal{P}_{e}(G)$ is $2K_2$-free;
\item $G$ is a cyclic group, a dihedral group, or an elementary abelian $2$-group.
\end{enumerate}
\end{thm}
\proof
It is easily seen that (d) implies (a) and (b). Since split graph and threshold graph are $2K_2$-free, any of (a) and (b) implies (c). Finally, by Lemma \ref{ep-sg-lem2}, (c) implies (d), which finishes the proof of the theorem.
\qed

%\begin{remark}
%\red{Example~\ref{co-dihe} implies that the enhanced power graph of a dihedral group is $P_4$-free.}
%\end{remark}

\section{Chordal graphs and cographs}\label{}

In this section, we investigate which groups have chordal enhanced power graph, and we study which groups have the enhanced power graph that is a cograph. In the first subsection, we classify all finite nilpotent groups whose enhanced power graph has the above mentioned properties, and in the second one, we turn our focus on non-nilpotent finite groups.

In \cite{Za}, it was proven that, if groups $G$ and $H$ have relatively prime orders, then $P_e(G\times H)$ is the strong graph product of $P_e(G)$ and $P_e(H)$; see \cite[Lemma 2.1]{Za}. Consequently, for $(g_1,h_1),(g_2,h_2)\in G\times H$, $(g_1,h_1)\neq(g_2,h_2)$, we have that $(g_1,h_1)\sime(g_2,h_2)$ if and only if both $\langle g_1,g_2\rangle$ and $\langle h_1,h_2\rangle$ are cyclic. Note that we may denote an element $(g,h)\in G\times H$ simply by $gh$. We shall use that fact in the proof of the following lemma.
%
%In this section, we classify all finite nilpotent groups whose enhanced  power graphs are chordal graphs and cographs.
%We also give some families of non-nilpotent groups whose enhanced power graphs are chordal graphs and cographs.

\begin{lemma}\label{ep-pgroup}
Let $\Gamma$ be a graph without any closed twins, let $G$ be a group, and let $n$ be relatively prime to $\lvert G\rvert$. Then $\mathcal P_e(G\times \mathbb Z_n)$ is $\Gamma$-free if and only if $\mathcal P_e(G)$ is $\Gamma$-free.
\end{lemma}

\proof
Since $\mathcal P_e(G)$ is an induced subgraph of $\mathcal P_e(G\times\mathbb Z_n)$, if $\mathcal P_e(G\times \mathbb Z_n)$ is $\Gamma$-free, then  $\mathcal P_e(G)$ is $\Gamma$-free as well. Therefore, it only remains to show the other implication of the lemma.

Suppose now that $X\subseteq G\times\mathbb Z_n$ induces $\Gamma$ in $\mathcal P_e(G\times\mathbb Z_n)$. Let $\pi_1:G\times\mathbb Z_n\rightarrow G$ denote the first projection, and let $H=\pi_1(X)$. Let us prove that $\pi_1$ embeds $\Gamma$ into $\mathcal P_e(G)$. Suppose the opposite, i.e. that $(h,k),(h,l)\in X$ for some $h\in H$ and $k,l\in \mathbb Z_n$, $k\neq l$. Let $(g,m)\in G\times\mathbb Z_n$. Then, by \cite[Lemma 2.1]{Za},
\begin{align*}
\langle hk,gm\rangle\text{ is cyclic }&\Leftrightarrow \langle h,g\rangle\text{ and }\langle k,m\rangle\text{ are cyclic}\\
&\Leftrightarrow \langle h,g\rangle\text{ and }\langle l,m\rangle\text{ are cyclic}\\
&\Leftrightarrow\langle hl,gm\rangle\text{ is cyclic.}
\end{align*}
Therefore, $hk$ and $hl$ have the same closed neighborhood in $\mathcal P_e(G\times\mathbb Z_n)$, which implies that they have the same closed neighborhood in $\Gamma$. However, since no pair of vertices of $\Gamma$ have same closed neighborhoods, it follows that $(h,k)=(h,l)$. Thus, we have obtained a contradiction, and therefore, we have shown that $\pi_1$ embeds $\Gamma$ into $\mathcal P_e(G)$. This completes the proof that $\mathcal P_e(G)$ contains $\Gamma$ as an induced graph.
\qed

\subsection{Nilpotent groups}

In this subsection, we focus on nilpotent groups. Before proving the main result, we need to show the following lemma first.

\begin{lemma}\label{ep-nil-group}
Let $G$ be a nilpotent group. Then the following are equivalent:
\begin{enumerate}
\item $G$ has at most one noncyclic Sylow subgroup;
\item $\mathcal{P}_{e}(G)$ is $P_4$-free;
\item $\mathcal{P}_{e}(G)$ is $C_4$-free.
\end{enumerate}
\end{lemma}

\proof
Let us show first that any of (b) and (c) implies (a). Let $P$ and $Q$ be two distinct noncyclic Sylow subgroups of $G$. Then $P\times Q$ is a subgroup of $G$. Note that the orders of $|P|$ and $|Q|$ are relatively prime. Since $P$ is not cyclic, neither is $\mathcal P_e(P)$ complete. Therefore, there are $x_1,x_2\in P$ that are nonadjacent in $\mathcal P_e(P)$. Similarly, there are $y_1,y_2\in Q$ nonadjacent in $\mathcal P_e(Q)$. Then, by \cite[Lemma 2.1]{Za}, $x_1\sime y_1\sime x_2\sime x_2y_2$ is an induced path of length four in $\mathcal P_e(G)$, and $x_1\sime y_1\sime x_2\sime y_2\sime x_1$ is an induced cycle of length four in $\mathcal P_e(G)$.

It remains now to show that (a) implies (b) and (c). Let us show first that the enhanced power graph of a finite $p$-group is $P_4$-free and $C_4$-free. Suppose that, for a $p$-group $P$, elements $x,y,z,w\in P$ induce the path $x\sime y\sime z\sime w$ or the cycle $x\sime y\sime z\sime w\sime x$  in $\mathcal P_e(P)$. Then $y\in\langle z\rangle$ or $z\in\langle y\rangle$, so let us assume, without loss of generality, that $y\in\langle z\rangle$. Then $z\in\langle x\rangle$ or $x,z\in\langle y\rangle$, which implies that $x\sime z$, which is a contradiction. Therefore, the enhanced power graph of any finite $p$-group is a cograph. Therefore, by Lemma \ref{ep-pgroup}, (a) does imply (b) and (c) since $P_4$ and $C_4$ contain no closed twins.
\qed

Now, we are ready to prove the main result of this subsection.

\begin{thm}\label{ep-cog-thm1}
The following are equivalent for any nilpotent group $G$:
\begin{enumerate}
\item $\mathcal{P}_{e}(G)$ is a chordal graph;
\item $\mathcal{P}_{e}(G)$ is a cograph;
\item $G$ has at most one noncyclic Sylow subgroup;
\item $G\cong P\times \mathbb{Z}_n$, where $P$ is a $p$-group for some $p$ relatively prime to $n$.
\end{enumerate}
\end{thm}

\proof
Clearly, conditions (c) are (d) are equivalent because every finite nilpotent group is the direct product of its Sylow subgroups.

Let us remind ourselves that a graph is a cograph if and only if it is $P_4$-free. Therefore, by Lemma \ref{ep-nil-group}, conditions (a), (b), and (c) are equivalent. This proves the theorem.
\qed

\begin{cor}\label{pgs}
The enhanced power graph of a finite $p$-group is a chordal graph, as well as a cograph.
\end{cor}

\begin{remark}
By Theorem \ref{ep-cog-thm1} and Lemma \ref{nnpg-lem-new}, if $G$ contains a nilpotent subgroup which has at least two distinct noncyclic Sylow subgroups, then $\mathcal{P}_{e}(G)$ is neither a chordal graph nor a cograph. For example, if $G$ has a subgroup isomorphic to $\mathbb{Z}_6\times \mathbb{Z}_6$, then $\mathcal{P}_{e}(G)$ is neither a chordal graph nor a cograph.
\end{remark}

Also, notice that some of the classes of graphs from Diagram \ref{hase-diag-1} are incomparable, such as cographs and chordal graphs. Interestingly however, by \cite[Theorem 6.2]{Za} and Theorems \ref{ep-thm-split} and \ref{ep-cog-thm1}, if we observe, instead of all graphs, only enhanced power graphs of finite nilpotent groups, we obtain the following relations between families of those graphs.
\begin{align*}
\text{threshold graphs}\subseteq\text{split graphs}\subseteq \text{chordal graphs}=\text{cographs}\subseteq\text{perfect graphs}
\end{align*}
%
% a different diagram in which any two classes of graphs are comparable. See Diagram \ref{hase-diag-2}.
%
%\begin{figure}[h!]
%\centering
%\begin{tikzpicture}
%  \node (a) at (0,0) {threshold graphs};
%  \node (b) at (0,1) {split graphs};
%  \node (c) at (0,2) {chordal graphs = cographs};
%  \node (d) at (0,3) {perfect graphs};
%  \draw (a) -- (b) -- (c) -- (d);
%\end{tikzpicture}
%\caption{}  \label{hase-diag-2}
%\end{figure}

\subsection{Non-nilpotent groups}

In this subsection, we give some families of non-nilpotent groups whose enhanced power graphs are chordal graphs and cographs (see Example \ref{co-dihe} and Propositions \ref{nng-p2}, \ref{nng-p3+nnng-p3}, \ref{An-co}, and \ref{An-chor}). First, we give a sufficient condition for a finite group whose enhanced power graph is chordal and a cograph. Let us remind the reader that $\mathcal M(G)$ denotes the set of all maximal cyclic subgroups of $G$.

\begin{prop}\label{nng-p1}
Let $G$ be a group. Let $k$ be a number such that the intersections of any two distinct maximal cyclic subgroups has cardinality $k$, then $\mathcal{P}_{e}(G)$ is a chordal graph as well as a cograph.
%Given a positive integer $k$, if $|M\cap N|=k$ for any two distinct $M,N\in \mathcal{M}_G$, then $\mathcal{P}_{e}(G)$ is a chordal graph as well as a cograph.
\end{prop}

\proof
Let $U=M\cap N$ for some $M,N\in\mathcal M(G)$. Since the intersection of two cyclic subgroups is a cyclic subgroup, $U$ is a cyclic subgroup $G$. Let $X$ and $Y$ be any maximal cyclic subgroups of $G$. We want to show now that $X\cap Y=U$. $X\cap M$ is also a cyclic subgroup of $M$ of order $k$. Because $M$ has the unique subgroup of order $K$, it follows that $U=X\cap M$. Similarly, as $X$ has the unique subgroup of order $k$, we obtain $X\cap Y=X\cap M=U$. Therefore, the intersection of any pair of distinct maximal cyclic subgroups of $G$ is $U$.

Suppose, for contradiction, that $x,y,y,w\in G$ induce a path $x\sime y \sime z \sime w$ or a cycle $x\sime y \sime z \sime w\sime x$ in $\mathcal{P}_{e}(G)$. Notice that elements of $U$ are adjacent to all elements of $G$ in $\mathcal P_e(G)$. Therefore, $U$ contains none of $x$, $y$, $z$, and $w$ since, otherwise, $\{x,y,z,w\}$ would not induce $P_4$ nor $C_4$. Let us denote two maximal cyclic subgroups that contain $\langle x,y\rangle$ and $\langle y,z\rangle$ by $K$ and $L$, respectively. Since $x\not\sime z$, it follows that $K\neq L$. However, $K\neq L$ implies that $y\in K\cap L\subseteq U$, which is a contradiction. Thus, $\mathcal P_e(G)$ is both $P_4$-free and $C_4$-free, and therefore, it is a chordal graph and a cograph.
%
%for any two $M,N\in \mathcal{M}_G$, where $U$ has size $k$.
%It follows that $\{x,y,z,w\}\cap U=\emptyset$.
%Since $\langle x,y\rangle$ is cyclic, there exists $K\in \mathcal{M}_G$ such that $\langle x,y\rangle\subseteq K$.
%Similarly, there exists $L\in \mathcal{M}_G$ such that $\langle y,z\rangle\subseteq L$. If $K\ne L$, then $y\in K\cap L$, and so
%$y\in U$, which is a contradiction. It follows that $K=L$, which implies that $\langle x,z \rangle$ is cyclic, this contradicts that  $x\sim y \sim z \sim w$ is an induced path. We conclude that $\mathcal{P}_{e}(G)$ is a cograph. Similarly, we can conclude that $\mathcal{P}_{e}(G)$ is $C_4$-free, and so  $\mathcal{P}_{e}(G)$ is chordal by Lemma \ref{obs1}.
\qed

%\medskip
%
%Let $p$ be a prime and \red{$n\geq 2$} a positive integer.
%Then the elementary abelian $p$-group $\mathbb{Z}_p^n$ has
%$\frac{p^n-1}{p-1}$ maximal cyclic subgroups isomorphic to $\mathbb{Z}_p$,
%and every two distinct maximal cyclic subgroups has trivial intersection. As a consequence, Proposition \ref{nng-p1} implies the following result
%(in fact, the following result also follows from Corollary \ref{pgs}).
%
%\begin{example}
%Let $p$ be a prime and let $n$ be a positive integer at least $2$. Then $\mathcal{P}_{e}(\mathbb{Z}_p^n)$
%is a chordal graph as well as a cograph.
%\end{example}

Notice that, for a group $G$, the condition from Proposition \ref{nng-p1} is equivalent to $G$ having a cyclic subgroup $C$ such that $C=M\cap N$ for every $M,N\in\mathcal M(G)$, $M\neq N$.

In a dihedral group, the intersection of any two distinct maximal cyclic subgroups is the trivial group. Also, a generalized quaternion group has the unique subgroup of order $2$, which is the intersection of every two distinct maximal cyclic subgroups. Therefore, by Proposition \ref{nng-p1}, we obtain the following which conclude that the enhanced power graph of any dihedral group or any generalized quaternion group is cograph and a chordal graph.

\begin{example}\label{co-dihe}
The enhanced power graph of every dihedral group and every generalized quaternion group is chordal and a cograph.
\end{example}

Let us remark that $D_{2n}$ and $Q_{4n}$ are nilpotent only if $n$ is a power of $2$.

A group is called a {\em $P$-group} \cite{De89} if every nonidentity element of the group has prime order. For example,
$D_{2q}$ is a $P$-group for any odd prime $q$. A group is called a {\em CP-group} \cite{Hi} if every nontrivial element of the group has prime power order. For example, any $p$-group is a $CP$-group, and $D_{2q^n}$ is a $CP$-group for any prime $q$ and positive integer $n$. Clearly, a $P$-group is also a $CP$-group. By Proposition \ref{nng-p1}, we see that for each $P$-group $G$, $\mathcal{P}_{e}(G)$ is a chordal graph as well as a cograph. Moreover, we have the following result.

\begin{prop}\label{nng-p2}
The enhanced power graph of any $CP$-group is a chordal graph as well as a cograph.
\end{prop}
\proof
Let $G$ be a $CP$-group. Suppose that some elements $x,y,z,w\in G$ induce the path $x\sime y \sime z \sime w$ or the cycle $x\sime y \sime z \sime w\sime x$ in $\mathcal{P}_{e}(G)$. Notice that $e\not\in\{x,y,z,w\}$ since the identity $e$ is adjacent to all vertices of $\mathcal P_e(G)$. Therefore, orders of $x$, $y$, $z$ and $w$ are powers of some prime $p$. Without loss of generality, we may assume that $z \in \langle y\rangle$. Then $x\sime z$ since either $x,z\in\langle y\rangle$ or $z\in\langle y\rangle\leq\langle x\rangle$. This contradicts the fact that $x$, $y$, $z$ and $w$ induce is a path or a cycle. Thus, $\mathcal{P}_{e}(G)$ is $P_4$-free and $C_4$-free, and therefore, $\mathcal P_e(G)$ is chordal and a cograph.
\qed

The {\em prime graph} of a finite group $G$ is a simple graph whose vertex set consists of all prime divisors of $|G|$ and whose two distinct vertices $p$ and $q$ are adjacent if $G$ has an element of order $pq$. The prime graph of a group was first introduced by Gruenberg and Kegel in an unpublished manuscript studying integral representations of groups in 1975. It is clear that, for a group $G$, its prime graph is a null graph if and only if $G$ is a $CP$-group. Thus, we have the following corollary.

\begin{cor}
Let $G$ be a group with null prime graph. Then $\mathcal{P}_{e}(G)$ is a chordal graph as well as a cograph.
\end{cor}

$CP$-groups were studied in \cite{Del,Suz,She}.
All nonabelian simple $CP$-groups were classified by
Suzuki \cite{Suz}. By \cite{Suz}, we have the following result.

\begin{example}
If $G$ is isomorphic to $Sz(8)$, $Sz(32)$, $PSL(3,4)$, or $PSL(2,q)$ for some $q\in\{5,7,8,9,17\}$, then $\mathcal{P}_{e}(G)$ is a chordal graph as well as a cograph.
\end{example}

In the remainder of this subsection, we classify all symmetric groups and all alternative groups whose enhanced power graph is a cograph and a chordal graph. Note that $S_n$ is non-nilpotent for every $n\ge 3$.
%
%
%The symmetric group of degree $n$, denoted by $S_n$,
%is the group of all permutations on $n$ symbols.
%Note that $S_n$ is non-nilpotent if $n\ge 3$.

\begin{prop}\label{nng-p3+nnng-p3}
The following conditions are equivalent.
\begin{enumerate}
\item $\mathcal P_e(S_n)$ is a cograph;
\item $\mathcal P_e(S_n)$ is a chordal graph;
\item $n\leq 5$.
\end{enumerate}
\end{prop}

\proof
Let us observe the subgraph of $\mathcal P_e(S_6)$ induced by $(1\ 2)$, $(3\ 4\ 5)$, $(1\ 6)$, and $(2\ 3\ 4)$. Notice that the subgroup generated by a pair of disjoint cyclic permutations of relatively prime orders is cyclic. On the other hand, if a pair cyclic permutations does not commute, they do not generate a cyclic subgroup. However, although a pair of disjoint cyclic permutations of the same order commute, they do not generate a cyclic subgroup either. Therefore,
\[(1\ 2)\sime(3\ 4\ 5)\sime(1\ 6)\sime(2\ 3\ 4)\]
is an induced path in $\mathcal P_e(S_6)$, and thus, $\mathcal P_e(S_6)$ is not a cograph. Therefore, by Lemma \ref{nnpg-lem-new}, (a) implies (c).

Similarly, notice that
\[\begin{gathered}
(1\ 2\ 3) \sime (5\ 6) \sime (2\ 3\ 4) \sime (6\ 1) \sime (3\ 4\ 5) \sime (1\ 2) \sime (4\ 5\ 6)\\
 \sime (2\ 3) \sime (5\ 6\ 1) \sime (3\ 4) \sime (6\ 1\ 2) \sime (4\ 5) \sime (1\ 2\ 3)
\end{gathered}\]
is an induced cycle of $\mathcal P_e(S_6)$. Therefore,
$\mathcal P_e(S_6)$ is not chordal either. Thus, by Lemma \ref{nnpg-lem-new}, (b) implies (c).

It remains for us to show that (c) implies (a) and (b). Suppose that $n\leq 5$, and suppose that $\mathcal P_e(S_n)$ contains elements $x$, $y$, $z$, and $w$ which induce the path $x\sime y\sime z\sime w$ or the cycle $x\sime y\sime z\sime w\sime x$. Notice that, in a path and in a cycle, no pair of vertices have the same closed neighborhood. Therefore, none of $x$, $y$, $z$, and $w$ is a cyclic permutation of order $5$ or $4$, and none of them is the product of two transpositions, since, for a cyclic permutation $c$ of order $5$ or $4$, $[c]_{\equive}=\langle c\rangle\setminus\{e\}$. Thus, $\{x,y,z,w\}$ contains only cyclic permutations of order $3$ and transpositions. Without loss of generality, let $o(y)=3$. But a cyclic permutation of order $3$ commutes with at most one transposition in $S_n$, $n\leq 5$. Thus, we have obtained a contradiction, and we conclude that $\mathcal P_e(S_n)$ is $P_4$-free and $C_4$-free for $n\leq 5$. This finishes our proof.
\qed

\begin{remark}
In the above proof, the length of the presented induced cycle in $\mathcal P_e(S_6)$ is $12$. However, $\mathcal P_e(S_6)$ does contain an induced cycle of length $6$. In particular,
\[\begin{gathered}
(1\ 2)(3\ 4)(5\ 6)\sime (1\ 4\ 5)(2\ 3\ 6) \sime (1\ 3)(2\ 5)(4\ 6) \sime (1\ 5\ 6)(2\ 4\ 3)\\
\sime (1\ 4)(2\ 6)(3\ 5)\sime (1\ 6\ 3)(2\ 5\ 4) \sime (1\ 2)(3\ 4)(5\ 6)
\end{gathered}\]
is one such induced cycle. Although that fact may be easier to verify for the induced cycle from the proof of Lemma \ref{nng-p3+nnng-p3} than for the above-mentioned one. Moreover, one can show that $6$ is the smallest length of an induced cycle in $\mathcal P_e(S_6)$.

Furthermore, for any $n\ge 7$, $\mathcal P_e(S_n)$ contains an induced cycle of length $4$, e.g.
\[(1\ 2\ 3)\sime (4\ 5) \sime (1\ 2\ 7) \sime (4\ 6)\sime (1\ 2\ 3).\]
\end{remark}

%\medskip

Now, we turn our attention to alternating groups. It is well known that $A_n$ is a simple group for each $n\ge 5$. Also, note that $A_n$ is non-nilpotent for every $n\ge 4$.

%The alternating group of degree $n$, denoted by $A_n$, is a group of even permutations on a set of length $n$. It is well known that $A_n$ is a simple group for each $n\ge 5$. Remark that $A_n$ is non-nilpotent for each $n\ge 4$.

\begin{prop}\label{An-co}
$\mathcal{P}_{e}(A_n)$ is a cograph if and only if $n\le 6$.
\end{prop}

\proof
Notice that
\[(1\ 2)(3\ 4)\sime(5\ 6\ 7)\sime(1\ 3)(2\ 4)\sime(1\ 2\ 3\ 4)(5\ 6)\]
is an induced path in $\mathcal P_e(A_7)$. Namely, $(1\ 3)(2\ 4)=\big( (1\ 2\ 3\ 4)(5\ 6) \big)^2$, and as discussed in the proof of Proposition \ref{nng-p3+nnng-p3}, disjoint permutations of relatively prime orders generate cyclic subgroups, and any pair of elements of a cyclic subgroup commutes. By Lemma \ref{nnpg-lem-new}, this proves one implication of the proposition.

In order to prove the other implication as well, notice that $A_6$ is a $CP$-group. Then, by Proposition \ref{nng-p2}, $\mathcal{P}_{e}(A_6)$ is a cograph. Therefore, $\mathcal P_e(A_n)$ is a cograph for every $n\leq 6$. This completes the proof.
\qed

\begin{prop}\label{An-chor}
$\mathcal{P}_{e}(A_n)$ is a chordal graph if and only if $n\le 7$.
\end{prop}

\proof
Notice that, for $n\geq 8$, $\mathcal{P}_{e}(A_n)$ contains the induced cycle
\[(1\ 2)(3\ 4)\sime (5\ 6\ 7) \sime (1\ 3)(2\ 4) \sime (5\ 6\ 8)\sime (1\ 2)(3\ 4)\]
of length four.

By Lemma \ref{nnpg-lem-new}, it now suffices to prove that $\mathcal{P}_{e}(A_7)$ has no induced cycles of length at least $4$. Suppose, for contradiction, that $\cdots \sime x\sime y \sime z \sime w\sime \cdots$ is an induced cycle of $\mathcal{P}_{e}(A_7)$, where $x,y,z,w$ are pairwise distinct. Notice that $\pi_e(A_7)=\{1,2,3,4,5,6,7\}$. It follows that each of $o(x),o(y),o(z),o(w)$ is $2$ or $3$. Without loss of generality, let $o(y)=2$ and let $y=(1\ 2)(3\ 4)$. Then $x,z\in\{(5\ 6\ 7),(5\ 7\ 6)\}\subseteq\langle(5\ 6\ 7)\rangle$, which implies that $x$ is adjacent to $z$. Thus, we have obtained a contradiction.
\qed

\subsection{Small groups}

In this subsection, we describe some classes of groups whose enhanced power graphs are chordal graphs and cographs. % (see Corollary \ref{smallg-cor1} and Lemmas \ref{smallg-pq} and \ref{smallg-two-primes}).
Further, we determine the only two groups of order at most $24$ whose enhanced power graph is not a cograph. We also show the minimal order of a group with non-chordal enhanced power graph is $36$.

The following result is an immediate consequence of Lemma \ref{ep-pgroup}, since neither $P_4$ nor a cycle of length at least $4$ has closed twins.

\begin{cor}\label{smallg-cor1}
Let $G$ be a group of order relatively prime to $n$. Then $\mathcal{P}_{e}(G\times\mathbb Z_n)$ is a cograph {\rm (}resp. a chordal graph{\rm )} if and only if $\mathcal{P}_{e}(G)$ is a cograph {\rm (}resp. a chordal graph{\rm )}.
\end{cor}

Notice that a group whose order is a product of two distinct primes is either a cyclic or a $CP$-group. Therefore, by Proposition \ref{nng-p2}, we have the following result.

\begin{prop}\label{smallg-pq}
If $|G|$ is a product of two distinct primes, then $\mathcal{P}_{e}(G)$ is a chordal graph as well as a cograph.
\end{prop}

Before heading to the next lemma, we note that the {\em cyclicizer} of a group $G$, denoted by $\Cyc(G)$, is defined by
\begin{equation*}\label{peq1}
\Cyc(G)=\{g\in G: \langle g,x \rangle \text{ is cyclic for any $x\in G$}\}.
\end{equation*}
Some authors, e.g. \cite{Ob92}, also call it {\em cycel} of a group.

\begin{prop}\label{smallg-two-primes}
Let $G$ be a group and let $p,q$ be distinct prime numbers. Then $\mathcal{P}_{e}(G)$ is a chordal graph as well as a cograph if one of the following holds:
\begin{enumerate}
\item $\pi_e(G)=\{1,p,q,pq\}$ and $G$ has the unique subgroup of order $p$ or $q$;
\item $\pi_e(G)=\{1,p,q,pq\}$ and either $G$ has a unique cyclic subgroup of order $pq$, or the intersection of all cyclic subgroups of order $pq$ has size $p$ or $q$;
\item $\pi_e(G)=\{1,p,q,pq,p^2\}$ and either $G$ has a unique cyclic subgroup of order $pq$, or the intersection of all cyclic subgroups of order $pq$ is $\langle a \rangle$ where $a\in \Cyc(G)$ and $o(a)\in \{p,q\}$.
\end{enumerate}
\end{prop}

\proof
In order to prove this lemma, let us suppose that $x$, $y$, $z$, and $w$ induce the path $x\sime y\sime z\sime w$ or the cycle $x\sime y\sime z\sime w\sime x$ in $\mathcal P_e(G)$. Further, we will suppose that (a), (b), and (c) hold, respectively, find a contradiction in each of the cases. Thus, we shall prove the lemma.

(a) Let $\pi_e(G)=\{1,p,q,pq\}$, and, without loss of generality, let $G$ have the unique subgroup of order $p$. It follows that $\{o(y),o(z)\}=\{p,q\}$. Without loss of generality, let $o(y)=p$ and $o(z)=q$. Then $o(w)\not\in \{p,q\}$, and therefore, $o(w)= pq$. Thus, $\langle y\rangle\leq\langle w\rangle$, since $\langle y\rangle$ is the unique cyclic subgroup of order $p$. Therefore, $y\sime w$, which is a contradiction.

(b) Let $\pi_e(G)=\{1,p,q,pq\}$, and suppose that $G$ has a unique cyclic subgroup of order $pq$. Then $y$ and $z$ have prime orders. Without loss of generality, let $o(y)=p$ and $o(z)=q$. It follows that $o(w)\in\{p,pq\}$. Since $G$ has a unique cyclic subgroup of order $pq$, then $\langle w,z\rangle=\langle y,z\rangle$ is the unique cyclic subgroup of order $pq$, which implies that $y\sime w$. Thus we have obtained a contradiction.

Now, let $\pi_e(G)=\{1,p,q,pq\}$, and let the intersection of all cyclic subgroups of order $pq$ has size $p$ or $q$. Again, $o(y),o(z)\in\{p,q\}$, and we can assume without loss of generality that $o(y)=p$ and $o(z)=q$. Since the intersection of all cyclic subgroups of order $pq$ has size $p$, then $y$ belongs to the cyclic subgroup $\langle w,z\rangle$ of order $pq$. Therefore, $y\sime w$, which is a contradiction.

(c) Suppose that $\pi_e(G)=\{1,p,q,pq,p^2\}$ and that $G$ has a unique cyclic subgroup of order $pq$. Since elements $y$ and $z$ have prime orders, we may, without loss of generality, assume that $o(y)=p$ and $o(z)=q$. Then $o(w)\in\{p,pq\}$. Note that $\langle y,z\rangle$ is a cyclic subgroup of order $pq$. Since $G$ has a unique cyclic subgroup of order $pq$, then clearly $y\sime w$, which is a contradiction.

Now, let $\pi_e(G)=\{1,p,q,pq,p^2\}$, and let, without loss of generality, the intersection of all cyclic subgroups of order $pq$ be $\langle a \rangle$ where $a\in \Cyc(G)$ and $o(a)=p$. Also, $y$ and $z$ have prime orders. Without loss of generality, let $o(y)=p$ and $o(z)=q$. Then $y\in\Cyc(G)$, which implies that $y\sime w$. Thus, we have obtained a contradiction. This proves the lemma.
\qed

\begin{example}
$\mathbb{Z}_2\times \mathbb{Z}_2 \times S_3$,  $\mathbb{Z}_3\times S_3$ and  $SL(2, 3)$ satisfy the conditions {\rm (a), (b)} and {\rm(c)} of Proposition~{\rm\ref{smallg-two-primes}}, respectively.
\end{example}

By GAP \cite{Gap}, we have determined that the minimal size of a group with non-chordal enhanced power graph is $36$. In particular, Theorem~\ref{ep-cog-thm1} implies that $\mathcal P_e(\mathbb Z_6\times\mathbb Z_6)$ is not chordal. Even though the enhanced power graph of a finite nilpotent group is chordal if and only if it is a cograph by Theorem \ref{ep-cog-thm1}, that is not the case for finite groups in general. Using GAP, we have determined that the minimal size of a group whose enhanced power graph is not a cograph is $24$. Namely, the two such groups of order $24$ are $Q_{12}\times\mathbb Z_2$ and the semidirect product $Q_8\rtimes\mathbb Z_3$ with action kernel $\mathbb{Z}_2\times \mathbb{Z}_2$.

\bigskip
\noindent \textbf{Acknowledgements}~~
The first author's research is supported by National Natural Science Foundation of China (Grant No. 11801441).
The second author acknowledges financial support of the Ministry of Education, Science and Technological Development of the Republic of Serbia (Grant No. 451-03-68/2022-14/200125).
The fourth author's research is supported by National Natural Science Foundation of China (Grant No. 61976244).

\end{CJK*}

\end{document}